\documentclass[11pt]{article}

\usepackage{mathptmx}       
\usepackage{helvet}         
\usepackage{courier}        
\usepackage{type1cm}        
\usepackage{makeidx}         
\usepackage{graphicx}        
\usepackage{multicol}        
\usepackage[bottom]{footmisc}

\usepackage{amsmath}
\usepackage{amssymb}
\usepackage{amscd}
\usepackage{indentfirst}

\newcounter{fig}
\newtheorem{prop}{Proposition}
\newcommand{\dem}{\noindent {\bf Proof\ }}
\newtheorem{defi}{Definition}
\newtheorem{ques}{Question}

\newcommand{\Vect}{{\rm Vect}}
\providecommand{\scal}[2]{\langle #1, #2 \rangle}

\newcommand{\fin}{\hfill$\Box$}
\newcommand{\fine}{\tag*{\mbox{$\Box$}}}

\newcommand{\virg}{\raisebox{.7mm}{,}}

\def\R{{\Bbb R}}
\def\N{{\Bbb N}}

\def\1{{\bf 1}}

\makeindex             

\begin{document}

\title{An arithmetical function related to\\ B\'aez-Duarte's criterion for the Riemann hypothesis}

\author{Michel Balazard}
\date{}

\maketitle

\begin{flushright}
\textit{To the memory of my friend, Luis B\'aez-Duarte.}
\end{flushright}

\abstract{\footnotesize In this mainly expository article, we revisit some formal aspects of B\'aez-Duarte's criterion for the Riemann hypothesis. In particular, starting from Weingartner's formulation of the criterion, we define an arithmetical function $\nu$, which is equal to the M\"obius function if, and only if the Riemann hypothesis is true. We record the basic properties of the Dirichlet series of $\nu$, and state a few questions.}

\begin{center}
  {\sc Keywords}
\end{center}
\begin{quote}
{\footnotesize Riemann hypothesis, arithmetical functions, Dirichlet series, Hilbert space }
\end{quote}


\section{The spaces $D$ and $D_0$}

We will denote by $\N$ (resp. $\N^*$) the set of non-negative (resp. positive) integers, by $H$ the Hilbert space $L^2(0,\infty ; t^{-2}dt)$, with inner product
\[
\scal{f}{g}=\int_0^{\infty}f(t)\overline{g(t)}\;\frac{dt}{t^2}\virg
\]
and by $\Vect({\mathcal F})$ the set of finite linear combinations of elements of a family ${\mathcal F}$ of elements of $H$.


For $k \in \N^*$, we define
\[
e_{k}(t)=\{t/k\} \quad (t>0),
\]
where $\{u\}=u-\lfloor u\rfloor$ denotes the fractional part of the real number $u$, and $\lfloor u\rfloor$ its integer part. The functions $e_{k}$ belong to $H$, as do the functions $\chi$ and $\kappa$ defined by
\[
\chi(t)=[t\ge 1] \quad ; \quad \kappa(t)=t[0 <t <1]
\]
(here, and in the following, we use Iverson's notation : $[P]=1$ if the assertion $P$ is true, $[P]=0$ if it is false).


Let $D$ be the closed subspace of functions $f \in H$ of the type 
\begin{equation}\label{t7}
f(t)=\lambda t+\varphi(t),
\end{equation}
where $\varphi$ is constant on each interval $[j,j+1[$, $j \in \N$ (for $j=0$, the constant must be $0$). The functions~$e_k$ belong to $D$.

Let $D_0$ be the subspace of $D$ defined by taking $\lambda=0$ in \eqref{t7}, that is, the subspace of functions~$\varphi \in~H$ which are constant on each interval $[j,j+1[$, $j \in \N$. The functions~$\chi$ and $e_k-e_1/k$ belong to $D_0$.

A hilbertian basis for $D_0$ is given by the family of step functions $\varepsilon_k$ defined by
\[
\varepsilon_k(t)=\sqrt{k(k+1)}\cdot [k\le t<k+1]\quad (k \in \N^*, \, t>0).
\]

The mapping $h \mapsto \big(h(j)\big)_{j\ge 1}$ is a Hilbert space isomorphism of $D_0$ onto the sequence space ${\mathfrak h}$ of complex sequences $(x_j)_{j\ge 1}$ such that
\[
\sum_{j \ge 1} \frac{\lvert x_j\rvert^2}{j(j+1)} < \infty \cdotp
\]

Observe that, for $f \in D$, written as \eqref{t7}, one has
\begin{align}
\lambda &=\scal{f}{\kappa}\label{180425c}\\
f&=\lambda e_1+h, \text{ where } h \in D_0.\label{180425d}
\end{align}

Thus, the subspace $D$ is the (non orthogonal) direct sum of $\Vect(e_1)$ and~$D_0$.

In formula \eqref{180425c}, the function $\kappa$ could be replaced by its orthogonal projection $\kappa'$ on $D$. The definition of the families $(\psi_n$) of Proposition \ref{180427a} and $(g_n)$ of Proposition \ref{180603a} below could be modified accordingly. We compute $\kappa'$ in the appendix.

\smallskip

To every function in $D$, one can associate certain arithmetical functions. Let~$f \in~D$, with $\lambda$ and~$h$ as in \eqref{180425c}, \eqref{180425d}. We first define the arithmetical function 
\begin{equation}\label{180423a}
u(n)=u(n;f)=-\lambda +h(n)-h(n-1) \quad ( n \in \N^*).
\end{equation}

With this definition, we see that the function $\varphi$ of \eqref{t7} is
given by
\[
\varphi(t)=-\lambda t+f(t)=-\lambda t+\lambda\{t\} +h(t)=\sum_{n\le t}u(n).
\]

Thus, $f(t)$ is the remainder term in the approximation of the sum function $\varphi(t)$ of the arithmetical function $u$ by the linear function $-\lambda t$. The fact that $f$ belongs to~$H$ implies, and is stronger than, the asymptotic relation $f(t)=o(t)$.

For $f \in D$, we will also consider the arithmetical function $w=\mu*u$, where $\mu$ denotes the M\"obius function,
\[
w(n)=w(n;f)=\sum_{d \mid n} \mu(n/d)u(d;f) \quad ( n \in \N^*).
\]

For instance,
\begin{equation*}
u(n;\chi)=[n=1] \quad ; \quad w(n;\chi)=\mu(n) \quad (n \in \N^*).
\end{equation*}

The arithmetical functions $u$ and $w$ depend linearly on $f$ and the correspondences are one-to-one.
\begin{prop}\label{180604d}
For $f \in D$,
\[
f=0 \Leftrightarrow u =0 \Leftrightarrow w =0.
\]
\end{prop}
\dem

The second equivalence follows from $w=u*\mu$ and $u=w*1$ (M\"obius inversion). It remains to prove that $u=0 \Rightarrow f=0$. By \eqref{180423a}, $u=0$ implies $h(n)=\lambda n$ for all $n$, hence $\lambda =0$ since $h \in D_0$, and~$h=0$.\fin

\smallskip

Since $u=w*1$, one has
\[
f(t)=\lambda t +\sum_{n\le t}u(n)=\lambda t +\sum_{n\ge 1}w(n)\lfloor t/n\rfloor.
\]

In Proposition \ref{180427b} below, we will prove the identity
\begin{equation}\label{180425b}
\sum_{n\ge 1} \frac{w(n)}n=-\lambda,
\end{equation}
so that, for every $f$ in $D$ and every $t>0$, one has
\begin{equation}\label{180604c}
f(t)=-\sum_{n \ge 1} w(n)e_n(t).
\end{equation}

Of course, it does not mean that the series $\sum_{n \ge 1} w(n;f)e_n$ converges in $H$ (in fact, it diverges if~$f=\chi$, cf. \cite{MR1692568}, Theorem 2.2, p. 6), but, if it does, its sum is $-f$.

\section{Vasyunin's biorthogonal system}

In Theorem 7 of his paper \cite{zbMATH00857404}, Vasyunin defined a family $(f_k)_{k\ge 2}$, which, together with the family~$(e_k-e_1/k)_{k \ge 2}$, yields a biorthogonal system in $D_0$, which means that
\begin{equation}\label{t3}
\scal{e_j-e_1/j}{f_k}=[j=k] \quad (j\ge 2, \, k\ge 2).
\end{equation}

We will recall Vasyunin's construction, which can be thought of as a Hilbert space formulation of M\"obius inversion, and add several comments.

\subsection{The sequence $(\varphi_k)$}

First one defines, for $k \in \N^*$, a step function $\varphi_k\in D_0$ by
\[
\varphi_k(t)=k(k-1)[k-1\le t<k]-k(k+1)[k\le t <k+1]
\]
(Vasyunin's $\varphi_k$ have the opposite sign, according to his definition for $e_k$). Thus
\[
\varphi_k= \sqrt{k(k-1)}\cdot \varepsilon_{k-1}-\sqrt{k(k+1)}\cdot\varepsilon_{k} \quad (k \in \N^*),
\]
with $\varepsilon_0=0$ by convention. One sees that the family $(\varphi_k)_{k \ge 1}$ is total in $D_0$.

One checks that
\begin{equation}\label{t8}
\scal{h}{\varphi_k}=h(k-1)-h(k) \quad  (k\in \N^*),
\end{equation}
for $h \in D_0$ with constant value $h(k)$ on $[k,k+1[$ ($h(0)=0$). In particular,
\[
\scal{e_j-e_1/j}{\varphi_k}=[j \mid k]-1/j \quad (j\ge 1, \, k\ge 1).
\]

Using the family $(\varphi_k)$, one can write the values $u(n;f)$, for $f \in D$, as scalar products.
\begin{prop}\label{180427a}
For $f \in D$, with $\lambda$ and $h$ as in \eqref{180425c}, \eqref{180425d}, one has
\[
u(n;f)=\scal{f}{\psi_n},
\]
where
\[
\psi_n=(\scal{e_1}{\varphi_n}-1)\kappa -\varphi_n \quad ( n \in \N^*).
\]

In particular, $f \mapsto u(n;f)$ is a continuous linear form on $D$, for every $n \in \N^*$. 
\end{prop}
\dem

By \eqref{180425c}, \eqref{180423a} and \eqref{t8}, one has
\begin{align*}
u(n;f)&=-\scal{f}{\kappa} - \scal{h}{\varphi_n} \\
&=-\scal{f}{\kappa} - \scal{f-\scal{f}{\kappa} e_1}{\varphi_n} \\
&=  -\scal{f}{\kappa} -\scal{f}{\varphi_n}+\scal{e_1}{\varphi_n}\scal{f}{\kappa}  \\
& =\scal{f}{\psi_n}         \quad ( n \in \N^*).\fine
\end{align*}

We compute the scalar product $\scal{e_1}{\varphi_n}$ in the appendix.

\smallskip

The next proposition describes the behavior of the series $\sum_k\varphi_k/k$.

\begin{prop}\label{t19}
The series 
\[
\sum_{k \ge 1} \frac{\varphi_k}k
\]
is weakly convergent in $D_0$, with weak sum $-\chi$.
\end{prop}
\dem

The partial sum
\[
\sum_{k \le K} \frac{\varphi_k}k
\]
is the step function with values
\begin{align*}
0 &\text{ on } (0,1) \text{ and } (K+1,\infty)\\
-1 &\text{ on } (1,K)\\
-(K+1)&\text{ on } (K,K+1)
\end{align*}

This partial sum is thus equal to $-\chi$ on every fixed bounded segment of $(0,\infty)$, if $K$ is large enough, and the norm of this partial sum in $H$ is the constant $\sqrt{2}$. The result follows.\fin

\subsection{The sequence $(f_k)$}

Vasyunin defined
\[
f_k=\sum_{d \mid k} \mu(k/d)\varphi_d \quad ( k \in \N^*).
\]
Equivalently,
\begin{equation*}
\varphi_k=\sum_{d \mid k} f_d \quad (k \in \N^* ),
\end{equation*}
by M\"obius inversion ; this implies that the family $(f_k)_{k \ge 1}$ is also total in $D_0$.

A slightly more general form of \eqref{t3}, namely 
\begin{equation}\label{t5}
\scal{e_j-e_1/j}{f_k}=[j=k]-[k=1]/j \quad (j,k \in \N^*),
\end{equation}
is proved by means of the identity
\[
\sum_{j \mid d \mid k}\mu(k/d)=[j=k].
\]

Using the family $(f_k)$, one can write the values $w(n;f)$, for $f \in D$, as scalar products.
\begin{prop}\label{180603a}
For $f \in D$, with $\lambda$ and $h$ as in \eqref{180425c}, \eqref{180425d}, one has
\[
w(n;f)=\scal{f}{g_n},
\]
where
\[
g_n=(\scal{e_1}{f_n}-[n=1])\kappa -f_n  \quad ( n \in \N^*).
\]

In particular, $f \mapsto w(n;f)$ is a continuous linear form on $D$, for every $n \in \N^*$. 
\end{prop}
\dem

By Proposition \ref{180427a}, one has
\begin{align*}
w(n;f)&=\sum_{d \mid n}\mu(n/d)u(d;f) \\
&=\scal{f}{\sum_{d \mid n} \mu(n/d)\psi_d} \quad (n \in \N^*). 
\end{align*}

Now,
\begin{align*}
\sum_{d \mid n} \mu(n/d)\psi_d &=\sum_{d \mid n} \mu(n/d)\big((\scal{e_1}{\varphi_d}-1)\kappa -\varphi_d\big)\\
&=(\scal{e_1}{f_n}-[n=1])\kappa -f_n.\fine
\end{align*}

We compute the scalar product $\scal{e_1}{f_n}$ in the appendix.

\smallskip

In order to study the series $\sum_kf_k/k$, we will need the following auxiliary proposition.

\begin{prop}\label{180423e}
Let
\[
f(x)=\sum_{k \le x} \eta(k) \quad (x > 0),
\]
where $\eta$ is a complex arithmetical function such that $\eta(k)=O(1/k)$, for $k \ge 1$.

Then, for every fixed $\alpha >1$, 
\[
\sum_{k\ge 1} \big\lvert f(x/k)-f\big(x/(k+1)\big)\big\rvert^{\alpha}=O(1) \quad (x > 0).
\]
\end{prop}
\dem

The series is in fact a finite sum, as 
\[
f(x/k)=f\big(x/(k+1)\big)=0 \quad (k>x).
\]

We will use the estimate
\[
f(y)-f(x) \ll \sum_{x<k\le y}\frac1k \,\ll\, \frac 1x + \ln (y/x) \quad (y>x\ge 1).
\]

Thus,
\[
f(x/k)-f\big(x/(k+1)\big) \ll \frac kx+\frac 1k \ll \frac 1k\quad (k\le \sqrt{x}),
\]
and
\[
\sum_{k \le \sqrt{x}}\big\lvert f(x/k)-f\big(x/(k+1)\big)\big\rvert^{\alpha} \ll \sum_{k \ge 1}\frac 1{k^{\alpha}} \ll 1 \quad (x >0).
\]

If $k > \sqrt{x}$, then 
\[
\frac xk-\frac x{k+1} <1,
\]
so that the interval $]x/(k+1),x/k]$ contains at most one integer, say $n$, and, if $n$ exists, one has~$k=~\lfloor x/n\rfloor$ and
\[
f(x/k)-f\big(x/(k+1)\big)=\eta(n)\ll \frac 1n\cdotp
\]

Hence
\[
\sum_{k > \sqrt{x}}\big\lvert f(x/k)-f\big(x/(k+1)\big)\big\rvert^{\alpha} \ll \sum_{n \ge 1}\frac 1{n^{\alpha}} \ll 1 \quad (x >0).
\]

The result follows.\fin

\begin{prop}\label{180423d}
The series 
\[
\sum_{k \ge 1} \frac{f_k}k
\]
is weakly convergent in $D_0$ (hence in $H$), with weak sum $0$.
\end{prop}
\dem

Let $K \in \N^*$. One has
\[
S_K=\sum_{k\le K} \frac{f_k}k=\sum_{d\le K} \frac{m(K/d)}d\varphi_d,
\]
where
\[
m(x)=\sum_{n\le x}\frac{\mu(n)}{n} \quad (x>0).
\]

Hence,
\begin{align*}
S_K &=\sum_{d\le K} \frac{m(K/d)}d\big(\sqrt{d(d-1)}\cdot \varepsilon_{d-1}-\sqrt{d(d+1)}\cdot\varepsilon_{d}\big)\\
&=\sum_{d\le K-1}\Big(\frac{m\big(K/(d+1)\big)}{d+1}-\frac{m(K/d)}d\Big)\sqrt{d(d+1)}\cdot\varepsilon_{d} -\sqrt{1+1/K}\cdot\varepsilon_{K}
\end{align*}

For every fixed $d \in \N^*$, the fact that $\scal{S_K}{\varepsilon_d}$ tends to $0$ when $K$ tends to infinity follows from this formula and the classical result of von Mangoldt, that $m(x)$ tends to $0$ when $x$ tends to infinity.

It remains to show that $\lVert S_K\rVert$ is bounded. One has
\begin{align*}
\lVert S_K\rVert^2 &=\sum_{d\le K-1}d(d+1)\Big(\frac{m(K/d)}d-\frac{m\big(K/(d+1)\big)}{d+1}\Big)^2 + 1+1/K\\
&\le 2\sum_{d\le K-1}d(d+1)\Big(\frac{m(K/d)-m\big(K/(d+1)\big)}d\Big)^2 \\
&\qquad+2\sum_{d\le K-1}d(d+1)\Big(\frac{m\big(K/(d+1)\big)}{d(d+1)}\Big)^2+ 1+1/K\\
&\ll 1+\sum_{d\le K-1}\Big(m(K/d)-m\big(K/(d+1)\big)\Big)^2 
\end{align*}

The boundedness of $\lVert S_K\rVert$ then follows from Proposition \ref{180423e}.\fin

\smallskip

We are now able to prove \eqref{180425b}.

\begin{prop}\label{180427b}
Let $f \in D$, with $\lambda$ and $h$ as in \eqref{180425c}, \eqref{180425d}. The series
\[
\sum_{n \ge 1} \frac{w(n;f)}n
\]
is convergent and has sum $-\lambda$.
\end{prop}
\dem

Putting $\beta_N= \sum_{n \le N}f_n/n$ for $N \in \N^*$, one has
\begin{align*}
\sum_{n \le N} \frac{g_n}n &=\sum_{n \le N} \frac{(\scal{e_1}{f_n}-[n=1])\kappa -f_n}n\\
&=\big(\scal{e_1}{\beta_N}-1\big)\kappa-\beta_N,
\end{align*}
which tends weakly to $-\kappa$, as $N$ tends to infinity, by Proposition \ref{180423d}.

Hence,
\begin{equation*}
\sum_{n \le N} \frac{w(n;f)}n =\sum_{n \le N} \frac{\scal{f}{g_n}}n
=\scal{f}{\sum_{n \le N}g_n/n} \rightarrow -\scal{f}{\kappa}=-\lambda \quad ( N\rightarrow \infty).\fine
\end{equation*}

\section{Dirichlet series} 

For $f \in D$ we define
\[
F(s)=\sum_{n\ge 1} \frac{u(n;f)}{n^s}\virg
\]
and we will say that $F$ is the Dirichlet series of $f$.

\smallskip

We will denote by $\sigma$ the real part of the complex variable~$s$. The following proposition summarizes the basic facts about the correspondance between elements $f$ of~$D$ and their Dirichlet series $F$. We keep the notations of \eqref{180425c} and \eqref{180425d}.

\begin{prop}\label{180604a}
For $f \in D$, the Dirichlet series $F(s)$ is absolutely convergent in the half-plane $\sigma >~3/2$, and convergent in the half-plane $\sigma >1$. It has a meromorphic continuation to the half-plane~$\sigma >1/2$ (we will denote it also by $F(s)$), with a unique pole in $s=1$, simple and with residue $-\lambda$. In the strip~$1/2< \sigma < 1$, one has 
\begin{equation}\label{t17}
F(s)/s=\int_0^{\infty}f(t)t^{-s-1} dt .
\end{equation}

If $f \in D_0$, that is $\lambda=0$, there is no pole at $s=1$, and the Mellin transform \eqref{t17} represents the analytic continuation of $F(s)/s$ to the half-plane $\sigma >1/2$. Moreover, the Dirichlet series $F(s)$ converges on the line $\sigma=1$.
\end{prop}
\dem

If $h=0$ in \eqref{180425d}, the arithmetical function $u$ is the constant $-\lambda$, and $F=-\lambda \zeta$. In this case, the assertion about \eqref{t17} follows from (2.1.5), p. 14 of \cite{titchmarsh-zeta}.

If $\lambda =0$, then $f=h \in D_0$ and $u(n)=h(n)-h(n-1)$ by \eqref{180423a}. Therefore,
\begin{align*}
\sum_{n\ge 1} \frac{\lvert u(n)\rvert}{n^{\sigma}} &\le 2\sum_{n\ge 1} \frac{\lvert h(n)\rvert}{n^{\sigma}} \\
&\le 2\big(\sum_{n\ge 1} \frac{\lvert h(n)\rvert^2}{n^{2}} \Big)^{1/2}\big(\sum_{n\ge 1} \frac{1}{n^{2\sigma -2}} \Big)^{1/2}\\
&\le 2\zeta(2\sigma-2)^{1/2} \| h \|<\infty,
\end{align*}
if $\sigma >3/2$, where we used Cauchy's inequality for sums.

The convergence of the series $F(1)$ follows from the formula $u(n)= -\scal{h}{\varphi_n}$ and Proposition \ref{t19}. It implies the convergence of $F(s)$ in the half-plane $\sigma >1$.

Using the Bunyakovsky-Schwarz inequality for integrals, and the fact that $h=0$ on~$(0,1)$, one sees that the integral \eqref{t17} now converges absolutely and uniformly in every half-plane $\sigma ~\ge~1/2+~\varepsilon$ (with $\varepsilon >0$), thus defining a holomorphic function in the half-plane $\sigma >1/2$. It is the analytic continuation of $F(s)/s$ since one has, for~$\sigma > 3/2$, 
\begin{align*}
\int_0^{\infty}h(t)t^{-s-1} dt  &=\frac 1s \sum_{n \ge 1} h(n)\big(n^{-s}-(n+1)^{-s}\big)\\
&=\frac 1s \sum_{n \ge 1} \frac{h(n)-h(n-1)}{n^{s}}=\frac{F(s)}s\cdotp
\end{align*}

Finally, the convergence of the Dirichlet series $F(s)$ on the line $\sigma =1$ follows from the convergence at $s=1$ and the holomorphy of $F$ on the line, by a theorem of Marcel Riesz (cf. \cite{zbMATH02611087}, Satz I, p.~350).

One combines the two cases, $h=0$ and $\lambda =0$, to obtain the statement of the proposition.\fin

\smallskip

The Dirichlet series $F(s)$ of functions in $D_0$ are precisely those which converge in some half-plane and have an analytic continuation to $\sigma >1/2$ such that $F(s)/s$ belongs to the Hardy space $H^2$ of this last half-plane. As we will not use this fact in the present paper, we omit its proof.

\smallskip

We now investigate the Dirichlet series 
\[
\frac{F(s)}{\zeta(s)}=\sum_{n\ge 1} \frac{w(n;f)}{n^s}\cdotp
\]

\begin{prop}\label{180605a}
Let $f \in D$, and let $F(s)$ be the Dirichlet series of $f$. The Dirichlet series $F(s)/\zeta(s)$ is absolutely convergent if $\sigma >3/2$, and convergent if $\sigma \ge 1$.
\end{prop}
\dem

The Dirichlet series $F(s)$ converges for $\sigma >1$, and converges absolutely for~$\sigma >3/2$ (Proposition~\ref{180604a}). The Dirichlet series $1/\zeta(s)$ converges absolutely for~$\sigma >1$. The Dirichlet product $F(s)/\zeta(s)$ thus converges absolutely for~$\sigma >~3/2$, and converges for $\sigma >1$. 

If $s=1$, the series is convergent by Proposition~\ref{180427b}. Since the function $F(s)/\zeta(s)$ is holomorphic in the closed half-plane $\sigma \ge 1$, Riesz' convergence theorem applies again to ensure convergence on the line $\sigma=1$.\fin

\section{B\'aez-Duarte's criterion for the Riemann hypothesis}\label{180604b}

We now define
\begin{equation*}
{\mathcal B}=\Vect(e_{n}, \, n \in \N^*)\quad ; \quad {\mathcal B}_0=\Vect(e_{n}-e_1/n, \, n \in \N^*, \, n \ge 2).
\end{equation*}

Since $e_n \in D$ and $e_n-e_1/n \in D_0$ for all $n \in \N^*$, one sees that 
\[
\overline{\mathcal B} \subset D \quad ; \quad \overline{{\mathcal B}_0} \subset D_0 \quad ; \quad \overline{{\mathcal B}_0} =\overline{\mathcal B}\, \cap D_0.
\]
The subspace $\overline{\mathcal B}$ is the (non orthogonal) direct sum of $\Vect(e_1)$ and~$\overline{{\mathcal B}_0}$.

\smallskip

We will consider the orthogonal projection $\widetilde{\chi}$ (resp. $\widetilde{\chi}_0$) of $\chi$ on $\overline{\mathcal B}$ (resp.~$\overline{\mathcal B_0}$). In 2003, B\'aez-Duarte gave the following criterion for the Riemann hypothesis.

\begin{prop}\label{t2}
The following seven assertions are equivalent.
\begin{align*}
 (i)&\quad \overline{{\mathcal B}} =D \quad ; \quad (i)_0 \quad \overline{{\mathcal B}_0} =D_0 \\
(ii) &\quad \chi \in \overline{{\mathcal B}} \quad ; \quad (ii)_0 \quad \chi \in \overline{{\mathcal B_0}} \\
 (iii)&\quad \widetilde{\chi}=\chi \quad ; \quad (iii)_0\quad \widetilde{\chi}_0=\chi\\
(iv) &\quad \text{the Riemann hypothesis is true.}\\
 \end{align*}
 \end{prop}

In fact, B\'aez-Duarte's paper \cite{MR2057270} contains the proof of the equivalence of $(ii)$ and~$(iv)$ ; the other equivalences are mere variations. The statements $(i)_0$, $(ii)_0$ and~$(iii)_0$ allow one to work in the sequence space ${\mathfrak h}$ instead of the function space~$H$; see \cite{zbMATH05184517} for an exposition in this setting.

\smallskip

The main property of Dirichlet series of elements of $\overline{{\mathcal B}}$ is given in the following proposition.

\begin{prop}\label{180605b}
If $f \in \overline{{\mathcal B}}$, the Dirichlet series $F(s)/\zeta(s)$ has a holomorphic continuation to the half-plane $\sigma >1/2$.
\end{prop}
\dem

Write $f=\lambda e_1 +h$, with $\lambda \in \R$ and $h \in D_0$. If $h=0$, one has~$F=~-~\lambda \zeta$ and the result is true. 

Now suppose $\lambda=0$. The function $h$ is the limit in $H$ of finite linear combinations, say $h_j$ ($j \ge 1$), of the $e_k-e_1/k$ ($k \ge 2$), when $j \rightarrow \infty$. The Dirichlet series of~$e_k-e_1/k$ is
\[
(k^{-1}-k^{-s})\zeta(s),
\]
so that the result is true for each $h_j$. It remains to see what happens when one passes to the limit. 

By the relation between the Dirichlet series of $h_j$ and the Mellin transform of~$h_j$, one sees that the Mellin transform of~$h_j$ must vanish at each zero~$\rho$ of $\zeta$ in the half-plane $\sigma >1/2$, with a multiplicity no less than the corresponding multiplicity of $\rho$ as a zero of $\zeta$. Thus
\begin{equation}\label{t20}
\int_1^{\infty}h_j(t)t^{-\rho-1}\ln^k t \,dt =0
\end{equation}
for every zero $\rho$ of the Riemann zeta function, such that $\Re \rho >1/2$, and for every non-negative integer~$k$ smaller than the multiplicity of $\rho$ as a zero of $\zeta$. When $j\rightarrow \infty$, one gets \eqref{t20} with $h_j$ replaced by $h$, which proves the result for $h$.

One combines the two cases, $h=0$ and $\lambda =0$, to obtain the statement of the proposition.\fin

\section{The $\nu$ function}

\subsection{Weingartner's form of B\'aez-Duarte's criterion}

For $N \in \N^*$, we will consider the orthogonal projections of $\chi$ on the subspaces $\Vect(e_1,\dots,e_N)$ and $\Vect(e_2-e_1/2,\dots,e_N-e_1/N)$ :
\begin{align}
\chi_N&=\sum_{k=1}^N c(k,N)e_k\label{t4}\\
\chi_{0,N}&=\sum_{k=2}^N c_0(k,N)(e_k-e_1/k),\label{u4}
\end{align}
thus defining the coefficients $c(k,N)$ and $c_0(k,N)$. 
In \cite{MR2304340}, Weingartner gave a formulation of B\'aez-Duarte's criterion in terms of the coefficients $c_0(k,N)$ of \eqref{u4}. The same can be done with the $c(k,N)$ of \eqref{t4}. First, we state a basic property of these coefficients.

\begin{prop}
For every $k \in \N^*$, the coefficients $c(k,N)$ in \eqref{t4} and $c_0(k,N)$ in~\eqref{u4} (here, with~$k\ge~2$) converge when $N$ tends to infinity.
\end{prop}
\dem

With the notations of \S\ref{180604b},
\begin{align*}
\widetilde{\chi} &=\lim_{N \rightarrow \infty}\chi_N\\
\widetilde{\chi}_0 &=\lim_{N \rightarrow \infty}\chi_{0,N},
\end{align*}
where the limits are taken in $H$.

Using the identity \eqref{180604c}, we observe that, for every $N \in \N^*$, 
\begin{align*}
c(k,N) &= -w(k;\chi_N) \quad \;\;\,(k \ge 1)\\
c_0(k,N) &= -w(k;\chi_{0,N}) \quad (k \ge 2),
\end{align*}

Therefore, Proposition \ref{180603a} yields, for every $k$,
\begin{align*}
c(k,N) &\rightarrow -w(k;\widetilde{\chi}) \quad \;\;\,(N \rightarrow \infty)\\
c_0(k,N) &\rightarrow -w(k;\widetilde{\chi_0}) \quad (N \rightarrow \infty).\fine
\end{align*}

\fin

\begin{defi}
The arithmetical functions $\nu$ and $\nu_0$ are defined by
\begin{align*}
\nu(n)&=w(n;\widetilde{\chi}) \\
\nu_0(n)&=w(n;\widetilde{\chi_0}).
\end{align*}
\end{defi}

Note that 
\[
\nu_0(1)=\lim_{N\rightarrow \infty}\sum_{2\le k \le N}\frac{c_0(k,N)}k=-\sum_{k \ge 2}\frac{\nu_0(k)}k\virg
\]
by Proposition \ref{180427b}.

\smallskip

We can now state B\'aez-Duarte's criterion in Weingartner's formulation.
\begin{prop}
The following assertions are equivalent.
\begin{align*}
 (i) &\quad \nu= \mu\\
 (ii) &\quad  \nu_0=\mu \text{ on } \N^*\setminus\{1\}\\
 (iii) &\quad \text{the Riemann hypothesis is true.}
 \end{align*}
\end{prop}
\dem

By B\'aez-Duarte's criterion, $(iii)$ is equivalent to $\chi=\widetilde{\chi}$. By Proposition \ref{180604d}, this is equivalent to~$w(n;\chi)=w(n;\widetilde{\chi})$ for all $n \ge 1$, that is, $\mu=\nu$.

Similarly, $(iii)$ implies $\mu=\nu_0$. Conversely, if $\mu(n)=\nu_0(n)$ for all $n \ge 2$, then~$w(n;\chi-~\widetilde{\chi_0})=0$ for $n \ge 2$, which means that $\chi-\widetilde{\chi_0}$ is a scalar multiple of~$e_1$. This implies $\chi=\widetilde{\chi_0}$ since $\chi$ and $\widetilde{\chi_0}$ belong to $D_0$.\fin
\fin

\subsection{The Dirichlet series $\sum_n\nu(n)n^{-s}$}

Since $\nu(n)=w(n;\widetilde{\chi})$, the following proposition is a corollary of Propositions \ref{180605a} and \ref{180605b}.
\begin{prop}
The Dirichlet series
\[
\sum_{n \ge 1} \frac{\nu(n)}{n^s}
\]
is absolutely convergent for $\sigma >3/2$, convergent for $\sigma \ge 1$, and has a holomorphic continuation to the half-plane $\sigma >1/2$.
\end{prop}

\section{Questions}

Here are three questions related to the preceding exposition.

\begin{ques}\label{t22}
Is it true that $\widetilde{\chi}=\widetilde{\chi}_0$?
\end{ques}

\begin{ques}
Let $f \in D$ such that the Dirichlet series $F(s)/\zeta(s)$ has a holomorphic continuation to the half-plane $\sigma >1/2$. Is it true that $f \in \overline{{\mathcal B}}\,$?
\end{ques}

A positive answer would be a discrete analogue of Bercovici's and Foias' Corollary 2.2, p. 63 of~\cite{MR768266}. 

\begin{ques}
Is the Dirichlet series
\[
\sum_{n \ge 1} \frac{\nu(n)}{n^s}
\]
  convergent in the half-plane $\sigma >1/2$?
\end{ques}

Another open problem is to obtain any quantitative estimate beyond the tautologies $\| \widetilde{\chi}-\widetilde{\chi}_N\| =~o(1)$ and  $\| \widetilde{\chi}_0-\widetilde{\chi}_{0,N}\| =o(1)$ ($N \rightarrow \infty$).

\section*{Appendix : some computations}

\subsection*{Scalar products}

{\bf 1.} One has
\begin{equation}\label{180602a}
\scal{e_1}{\varepsilon_k}=\sqrt{k(k+1)}\int_k^{k+1}(t-k) \frac{dt}{t^2}
=\sqrt{k(k+1)} \big(\ln (1+1/k)-1/(k+1)\big).
\end{equation}

{\bf 2.} For $k \in \N^*$, one has
\begin{align*}
\scal{e_1}{\varphi_k} &=\int_{k-1}^k k(k-1)(t-k+1) \frac{dt}{t^2}-\int_{k}^{k+1} k(k+1)(t-k) \frac{dt}{t^2}\\
&=2k^2\ln k-k(k-1)\ln(k-1)-k(k+1)\ln(k+1) +1\\
&=-\omega(1/k),
\end{align*}
where
\begin{align*}
\omega(z) &=z^{-2}\big((1-z)\ln(1-z)+(1+z)\ln(1+z)\big)-1\\
&=\sum_{j \ge 1} \frac{z^{2j}}{(j+1)(2j+1)} \quad ( \lvert z \rvert \le 1).
\end{align*}

{\bf 3.} For $n \in \N^*$, one has
\begin{align*}
\scal{e_1}{f_n} &=\sum_{k \mid n}\mu(n/k)\scal{e_1}{\varphi_k}
=-\sum_{k \mid n}\mu(n/k)\omega(1/k)\\
&=-\sum_{j \ge 1} \frac{\sum_{k \mid n}\mu(n/k)k^{-2j}}{(j+1)(2j+1)}
=-\sum_{j \ge 1} \frac{n^{-2j}\prod_{p \mid n}(1-p^{2j})}{(j+1)(2j+1)}\cdotp 
\end{align*}

In particular,
\[
\sup_{n \in \N^*} \lvert \scal{e_1}{f_n} \rvert =\sum_{j \ge 1} \frac{1}{(j+1)(2j+1)}=\ln 4-1.
\]

\subsection*{Projections}

By \eqref{180602a}, the orthogonal projection $e'_1$ of $e_1$ on $D_0$ is
\[
e'_1=\sum_{k\ge 1} \scal{e_1}{\varepsilon_k}\varepsilon_k=\sum_{k\ge 1} \sqrt{k(k+1)}\big(\ln (1+1/k)-1/(k+1)\big)\varepsilon_k.
\]

Since $e'_1(k)$ has limit $1/2$ when $k$ tends to infinity, one sees that $e_1-e'_1$ "interpolates" between the fractional part (on $[0,1[\,$) and the first Bernoulli function (at infinity). One has the hilbertian decomposition
\[
D=D_0\oplus \Vect(e_1-e'_1).
\]

Since $\kappa \,\bot\, D_0$ and $\scal{\kappa}{e_1}=1$, the orthogonal projection of $\kappa$ on $D$ is
\[
\kappa'=\frac{e_1-e'_1}{\lVert e_1-e'_1\rVert^2}\cdotp
\]

\begin{center}
{\sc Acknowledgements}
\end{center}
{\footnotesize I thank Andreas Weingartner for useful remarks on the manuscript.}

\providecommand{\bysame}{\leavevmode ---\ }
\providecommand{\og}{``}
\providecommand{\fg}{''}
\providecommand{\smfandname}{et}
\providecommand{\smfedsname}{\'eds.}
\providecommand{\smfedname}{\'ed.}
\providecommand{\smfmastersthesisname}{M\'emoire}
\providecommand{\smfphdthesisname}{Th\`ese}

\medskip

\footnotesize

\noindent BALAZARD, Michel\\
Aix Marseille Univ, CNRS, Centrale Marseille, I2M, Marseille, France\\
Adresse \'electronique : \texttt{balazard@math.cnrs.fr}


\begin{thebibliography}{1}

\bibitem{MR1692568}
{\scshape L.~B{\'a}ez-Duarte} -- {\og A class of invariant unitary
  operators\fg}, \emph{Adv. Math.} \textbf{144} (1999), p.~1--12.

\bibitem{MR2057270}
\bysame , {\og A strengthening of the {N}yman-{B}eurling criterion for the
  {R}iemann hypothesis\fg}, \emph{Atti Accad. Naz. Lincei Cl. Sci. Fis. Mat.
  Natur. Rend. Lincei (9) Mat. Appl.} \textbf{14} (2003), p.~5--11.

\bibitem{zbMATH05184517}
{\scshape B.~{Bagchi}} -- {\og {On Nyman, Beurling and Baez-Duarte's Hilbert
  space reformulation of the Riemann hypothesis.}\fg}, \emph{{Proc. Indian
  Acad. Sci., Math. Sci.}} \textbf{116} (2006), p.~139--146.

\bibitem{MR768266}
{\scshape H.~Bercovici {\normalfont \smfandname} C.~Foias} -- {\og A real
  variable restatement of {R}iemann's hypothesis\fg}, \emph{Israel J. Math.}
  \textbf{48} (1984), no.~1, p.~57--68.

\bibitem{zbMATH02611087}
{\scshape M.~{Riesz}} -- {\og {Ein Konvergenzsatz f\"ur {\it Dirichlet}sche
  Reihen.}\fg}, \emph{{Acta Math.}} \textbf{40} (1916), p.~349--361.

\bibitem{titchmarsh-zeta}
{\scshape E.~C. Titchmarsh} -- \emph{The theory of the {R}iemann zeta
  function}, second \smfedname, Clarendon Press, Oxford, 1986, revised by {D.
  R. H}eath-{B}rown.

\bibitem{zbMATH00857404}
{\scshape V.~{Vasyunin}} -- {\og {On a biorthogonal system associated with the
  Riemann hypothesis.}\fg}, \emph{{St. Petersbg. Math. J.}} \textbf{7} (1996),
  p.~405--419.

\bibitem{MR2304340}
{\scshape A.~Weingartner} -- {\og On a question of {B}alazard and {S}aias
  related to the {R}iemann hypothesis\fg}, \emph{Adv. Math.} \textbf{208}
  (2007), p.~905--908.

\end{thebibliography}
\end{document}